\documentclass[12pt]{article}
 \usepackage{color}

\usepackage{latexsym}
\usepackage{epic,eepic,epsfig}
\usepackage{amssymb}
\usepackage{amsfonts,srcltx}

\newtheorem{theorem}{Theorem}[section]
\newtheorem{prop}[theorem]{Proposition}
\newtheorem{lem}[theorem]{Lemma}
\newtheorem{rem}[theorem]{Remark}

\newtheorem{defi}[theorem]{Definition}
\newtheorem{exam}[theorem]{Example}

\newcommand{\mod}{\, \hbox{mod} \, }

\newcommand{\Aut}{\hbox{Aut}\, }

\newcommand{\'}{\' \i}

\newcommand{\qed}{\hfill \mbox{\raisebox{0.7ex}{\fbox{}}}}

\def\demo{{\bf Proof}\hskip10pt}
  \def\a{\alpha}
\def\s{\sigma}

 \def\MM{{\cal M}}

\def\og{\overline G} \def\oh{\overline H} \def\oc{\overline C}
 \def\oq{\overline Q}

 \def\oa{\overline A}
    \def\ol{\overline L} 
 \def\op{\overline P} \def\os{\overline S}
\def\ot{\overline T}  

\def\o1{\overline 1} \def\olh{\overline h}

\def\o{\overline}   \def\olr{\overline r} \def\oll{\overline \ell} \def\olt{\overline t }

\def\di{\bigm|} \def\lg{\langle} \def\rg{\rangle}
\def\nd{\mathrel{\bigm|\kern-.7em/}}

\def\f{\noindent}

\def\PSL{\hbox{\rm PSL}}
\def\PSU{\hbox{\rm PSU}}
 
\def\Aut{\hbox{\rm Aut}}
\def\Inn{\hbox{\rm Inn}}
\def\Syl{\hbox{\rm Syl}}

\def\Mon{\hbox{\rm Mon}}
\def\Aut{\hbox{\rm Aut}}
\def\mod{\hbox{\rm mod }}

\def\PGammaL{\hbox{\rm P$\Gamma $L}}
\def\GL{\hbox{\rm GL}} \def\PG{\hbox{\rm PG}}
\def\char{\hbox{\rm \,char\,}}
\newcommand{\B}[1]{\mathbb #1}

\begin{document}

\baselineskip=16pt
\newcounter{teocounter}
\addtocounter{teocounter}{1}
\newcounter{defcounter}
\addtocounter{defcounter}{1}
\newcounter{figcounter}
\addtocounter{figcounter}{1}
\newcounter{itc}
\addtocounter{itc}{1}
\newcounter{Lorrms}
\setcounter{Lorrms}{2} \vspace{2.0 cm}

\begin{center}

{\Large   Orientably-Regular $p$-Maps  and  Regular $p$-Maps}
\end{center}

\begin{center}  Shaofei Du, Yao Tian$^*$ and Xiaogang Li
\end{center}

\begin{center}
{\footnotesize  School of Mathematical Sciences\\ Capital Normal University \\Beijing, 100048, P.R.China}
\end{center}

\begin{abstract}
A map is called a {\it $p$-map} if it has a  prime $p$-power vertices.
 An  orientably-regular (resp. A regular )  $p$-map  is called {\it solvable} if the group $G^+$ of all orientation-preserving automorphisms (resp. the group $G$ of  automorphisms)
   is solvable; and called {\it normal} if  $G^+$ (resp.  $G$) contains the normal  Sylow $p$-subgroup.

 In this paper, it will be proved that  both   orientably-regular $p$-maps  and  regular   $p$-maps are  solvable and  except for  few cases that $p\in \{2, 3\}$, they are    normal.
 Moreover, nonnormal   $p$-maps  will be  characterized  and   some   properties   and constructions of   normal   $p$-maps  will be given.
\end{abstract}


\renewcommand{\thefootnote}{\empty}
\footnotetext{$^*$corresponding author}
\footnotetext{ Email:  tianyao202108@163.com }

\section{Introduction}
A {\em map} is a cellular decomposition of a closed surface.
An alternative way to describe  maps is to consider them as cellular embeddings of graphs into closed surfaces.
By an {\em automorphism of a map} ${\cal M}$ we mean an  automorphism of the {\em underlying graph} ${\cal G}$
which extends to a   self-homeomorphism of the surface.
These automorphisms  form a subgroup $\Aut({\cal M})$ of the automorphism  group $\Aut({\cal G})$ of $\cal G$.
It is well-known that  $\Aut({\cal M})$ acts semi-regularly on the set of all flags (in most cases, which are incident vertex-edge-face triples).
If the action is regular, then we call the map as well as the corresponding embedding  {\em  regular}.

In the case of orientable supporting surface, if the group $\Aut^+(\MM)$ of all orientation-preserving
automorphisms of $\MM$ acts regularly on the set of arcs (incident vertex-edge pairs) of $\MM$, then  we call $\MM$ an {\emph{orientably-regular\/} map.
Such maps fall into two classes: those that admit also orientation-reversing automorphisms,
which are called \emph{reflexible}, and those that do not, which are  called \emph{chiral}. Therefore, a reflexible map is a regular map but a chiral map is not.

One of the central problems in topological graph theory is to classify all the regular  maps of  given underlying graphs.
In a general setting, the classification problem was treated in \cite{GNSS}.
However, the problem has been solved only for a few particular classes of graphs.
For instance,  complete graphs,  complete bipartite graphs, complete multipartite graphs,  $n$-dimensional cubes $Q_n$ and so on.
Moreover, there are closed connections between  regular maps and  other mathematical theories such as group theory, hyperbolic geometry and complex curves and so on.
During the past  forty  years, there have been abundant references on regular maps, see  \cite{JaJo,Jon2,Wil} for an overview.
\begin{defi}
A map is called a {\it $p$-map} if the number of vertices is  $p^k$, where $p$ is   prime and $k\ge 1$.
\end{defi}

  In this paper, we shall concentrate on   orientably-regular $p$-maps and  regular   $p$-maps.

Many known orientably-regular maps and regular maps are  $p$-maps. For instance, every orientably-regular map of a complete graph must be  a $p$-map, see \cite{JaJo}.
 Every $n$-dimensional hypercube $Q_{n}$ contains $2^n$ vertices,  whose  orientably-regular maps and nonorientable regular maps  were classified in \cite{CCDKNW} and \cite{KN1}, respectively.
 Every regular map of complete bipartite graph $K_{2^e,2^e}$ has $2^{e+1}$ vertices and a classification of  orientably-regular embeddings of this graph played an important role in a final classification  of $K_{n,n}$,  see \cite{DJKNS1,DJKNS2,Jon1,KK3}. Orientably-regular maps whose automorphism group $G$ is  nilpotent can be concluded to that when $G$ is a 2-group  and of course, these maps have a 2-power vertices, see \cite{CDNS,MNS}.

 One will see that the studies of {\it $p$-maps} is very closed to that of finite $p$-groups and its automorphisms as well. In other words, the studies of $p$-maps  stimulate us to focus on  some related $p$-groups and
 it might   pose some new research problems  for group theorists.

 Regular embeddings of simple graphs of order $p$, $p^2$ and $p^3$ have been classified, see \cite{DKN1,DK2,ZDM1}. During the studying of these maps,  we obtained some new observations, ideas, methods and conjectures.
 this motivates us to pay more attention on   $p$-maps. The aim of this paper is to give  a basic theoretical  characterization  for  orientably-regular $p$-maps and regular   $p$-maps.

  As usual, by  $\B{D}_k$,  $O_p(G)$ and $H\lhd G$,  we denote the dihedral group of order $k$, the maximal normal $p$-subgroup of $G$, and   a normal subgroup $H$ of  $G$,  respectively.
 First  we introduce   the following   concepts.

\begin{defi}  An  orientably-regular (resp. A regular )  $p$-map $\MM $ is called {\it solvable}  if  $\Aut^+(\MM)$  (resp. $\Aut(\MM)$) is solvable; and called {\it normal} if
 $\Aut^+(\MM)$ (resp.  $\Aut(\MM))$ contains the normal  Sylow $p$-subgroup.
\end{defi}

\begin{rem}
  Suppose that $\MM$ is  an orientably-regular $p$-map which is reflexible. Then it is regular. Since $|\Aut(\MM):\Aut^+(\MM)|=2$, it follows that   $\Aut(\MM)$ is solvable if and only if
  $\Aut^+(\MM)$ is solvable. Take  a Sylow $p$-subgroup $P$ of $\Aut(\MM)$.  If  $p\ne 2$, then  $P\lhd  \Aut(\MM)$ if and only if $P\lhd \Aut^+(\MM)$.
   If $p=2$,  then $P\lhd \Aut(\MM)$ implies that $P\cap \Aut^+(\MM)$ is the normal Sylow $2$-subgroup of   $\Aut^+(\MM)$. However,  the converse is not true.
       For   example, consider the regular tetrahedron: the Sylow $2$-subgroup of $\Aut^+(\MM)\cong A_4$ is normal   but every Sylow $2$-subgroup of  $\Aut(\MM)\cong S_4$ is nonormal.
      \end{rem}

Now we are ready to state the main theorem of this paper.
\begin{theorem}\label{main}  Let $\MM$ be an orientably-regular $p$-map  or a  regular   $p$-map. Then
 \begin{enumerate}
   \item[{\rm (1)}] $\MM$ is solvable;
     \item[{\rm (2)}] $\MM$ is normal,  except  for the following two cases:
     \begin{enumerate}
     \item[{\rm (2.1)}] $p=2$,  $G/O_2(G)\cong \B{Z}_m\rtimes \B{Z}_2$ or  $\B{Z}_m\rtimes \B{D}_4$, where $m\ge 3$ is odd.
     \item[{\rm (2.2)}]  $p=3$,    $G/O_3(G)\cong S_4$.
     \end{enumerate} \end{enumerate}
\end{theorem}

It follows from Theorem~\ref{main} that   an orientably-regular  $p$-map  or a  regular   $p$-map can be  nonnormal  only if $p=2$ or 3. To give a description in details for these  nonnormal $p$-maps, we
introduce three  families of  (multi)graphs and  maps as well.

   (i) By ${\cal D}_d$  we denote the $d$-dipoles,  which   is the graph of two vertices  joined by $d$ parallel edges. Every orientably-regular map of ${\cal D}_d$ is isomorphic to ${\cal D}(d, e):=\MM(G; x, x^ey)$ (see \cite{NS}), where $G=\lg x,y\di x^d=y^2=1, x^y=x^e\rg$ and $e^2\equiv 1(\mod m).$

   (ii)  By ${\cal S}_d$  we denote the $d$-semistars,   which   is the single vertex with $d$ semi-edges. By $DM(d)$ and $EM(d)$, we denote the regular map of  ${\cal S}_d$ in a disc and a sphere, respectively (see \cite{LS}). Note that in this paper  we only consider surfaces without boundary, except for the case in here we include the disc as a surface).

   (iii) By   $C_n^{(m)}$  we denote the graph resulting from the cycle $C_n$ of length $n$ by replacing each
edge with $m$ parallel edges.  There exists a unique nonorientable regular map  of $C_3^{(2)}$, denoted by  ${\cal C}(3,2)$ (see \cite{HNSW}).

\begin{theorem}\label{main1}  Suppose that $\MM$ is  a nonnormal  orientably-regular $p$-map  or a  nonnormal  regular   $p$-map. Let $G=\Aut^+(\MM)$ or $\Aut(\MM)$.
Then the quotient   map ${\overline \MM}$ induced by $O_p(G)$  is one of the following maps:
 \begin{enumerate}
     \item[{\rm (1)}] $p=2$,   ${\overline \MM}={\cal D}(m,e)$, where $m\ge 3$ is odd and $e^2\equiv 1(\mod m)$ but $e\not\equiv 1(\mod m)$.
     $\MM$ is   a nonnormal orientably-regular 2-map, more precisely, it is either  chiral   or   reflexible and nonnormal  regular;
   \item[{\rm (2)}] $p=2$,   ${\overline \MM}=DM(m)$ and  $\MM$ is nonorientable and nonnormal regular;
    \item[{\rm (3)}] $p=2$,   ${\overline \MM}=EM(m)$ and  $\MM$ is normal  orientably-regular  but nonnormal regular;
     \item[{\rm (4)}]  $p=3$,    ${\overline \MM}={\cal C}(3,2)$ and $\MM$ is   nonorientable and nonnormal  regular.
     \end{enumerate}
\end{theorem}

Theorems~\ref{main} and ~\ref{main1} will be proved by combining several lemmas in Sections 3-6.
\begin{rem}
  For a normal orientably-regular $p$-map,   $G=P\!\rtimes \B{Z}_{m}$; and   for a normal regular $p$-map,  either $G=P\!\rtimes \B{D}_{2m}$  if $p$ is odd or   $G=P\!\rtimes \B{Z}_{m}$ if $p=2$.
In all cases, $m$ and $p$ are co-prime.  Therefore, for normal $p$-maps,   the studies   of  $p$-maps are essentially determinations of  automorphisms of a given $p$-group. Of course, it is  very difficult, because of  complexities of finite $p$-groups. In some senses, Theorem~\ref{main} is a  basic theoretical characterization of  orientably-regular $p$-maps and  regular $p$-maps, and it could  be a
starting point for studying such maps.
  \end{rem}

After this introductory section,  a brief description for regular maps and some   known results    will be given in Section 2; the solvability of  orientably-regular $p$-maps
and regular $p$-maps will be proved in Section 3; the normality of such maps will be discussed in Sections 4-6 according to $p\ge 5$, $p=3$ and $p=2$, respectively;  such $p$-maps that the induced group $G^V$ of $G$ on
   the set $V$ of vertices acts primitively will be dealt with in Section 7.  In addition, some examples of $p$-maps are presented  in the related sections.

Except for the notations mentioned as above,  more notations and terminologies     used in this paper are listed below.

\vskip 3mm

$|G|$, $|g|$:   the order of a group $G$ and the order of an element  $g$ in $G$, respectively;

 $|G:H|$,   $H\lhd G$ and $H\char G$:  the index of $H$ in $G$,    a normal subgroup  $H$ of $G$ and  a  characteristic subgroup $H$ of $G$, respectively;

$H_G$: the maximal normal subgroup of $G$ contained in the subgroup $H$;

$K\rtimes H$:  the semi-product of $K$ by $H$ where $K$ is normal;

$[a,b]$ and $[H,K]$: the commutator of two elements $a$ and $b$, and the  subgroup generated by all commutators $[h,k]$ where $h\in H$ and $k\in K$,  respectively;

$\pi '$: the complementary set of a subset   $\pi $ of   $\mathbb{P}$,  where $\mathbb{P}$ is the set of primes;

$\pi$-subgroup $H$:  a subgroup $H$ such that  every prime divisor of $|H|$  is contained in a subset $\pi $ of  $\mathbb{P}$;

$\pi$-Hall subgroup $H$: a $\pi$-subgroup $H$ such that   $(|H|, |G|/|H|)=1$;

$\Syl_p(G):$  the set of Sylow $p$-subgroups of $G$;

$\Phi(G)$:  the Frattini subgroup of $G$ (the intersection of all maximal subgroups of $G$);

$F(G)$: the Fitting subgroup of $G$ (the product of all  nilpotent normal  subgroups of $G$);

 $O_{\pi}(G)$: the maximal normal $\pi$-subgroup of $G$ (the intersection of maximal $\pi$-subgroups of $G$);

$\rm exp(G)$: the exponent of $G$, which is defined to be the least positive integer $s$ such that $x^s=1$ for all $x\in G$;

$\Inn(g)$: the inner automorphism of $G$ induced by $g\in G$.

\section{Preliminary Results}
In this section,  we shall   give a brief description  for  orientably-regualr and  regular maps and list some known  group theoretical   results used in this paper.

 \subsection{Regular maps  and orientably-regular maps}
(1) {\bf Regular Maps.}
\vskip 3mm
 A regular map can be described in the following way.
\begin{defi}\label{map1}
{\rm For a given finite set $F$ and three fixed-point-free involutory permutations $t, r, \ell $ on $F$, a quadruple
$\MM=\MM(F; t, r, \ell )$ is called a {\it combinatorial map} if they satisfy two conditions:
(1)\ $t\ell =\ell t$; (2)\ the group $\lg t,r, \ell \rg $ acts transitively on $F.$}
\end{defi}

For a given combinatorial map $\MM=\MM(F; t, r, \ell ),$  $F$  is called the {\it flag} set,
$t,  r, \ell$ are called {\it transversal, rotary, and longitudinal involution,} respectively.
The group $\lg t, r, \ell \rg $ is called the {\it monodromy group} of $\MM$, denoted by $\Mon(\MM)$.
We define the {\it vertices, edges} and {\it face-boundaries} of $\MM$ to be the
orbits of the subgroups $\lg t, r\rg$, $\lg t, \ell \rg $ and $\lg r, \ell \rg $, respectively.
The incidence in $\MM $ is defined by a nontrivial set  intersection.
To ensure that the underlying graph has no multiple edges, we let $\lg t, r\rg \cap \lg t, r\rg ^{\ell}=\lg t\rg $.

The map $\MM$ in Definition~\ref{map1} is {\it unoriented}. Clearly, the even-word subgroup $\lg tr, r\ell \rg $
of $\Mon(\MM)$ has  index at most 2.  If the index is 1, then $\MM$ is said to be {\it nonorientable}.
If the index is 2, then one may fix an orientation for $\MM $ and so $\MM $ is said to be
{\it orientable}, while  the group $\Aut^+(\MM)$ of all orientation-preserving automorphisms
  acts transitively  on  all arcs is exactly $\lg tr, r\ell\rg .$

Given two maps $\MM _1=\MM(F_1; t_1, r_1, \ell _1)$ and $\MM_2=\MM_2(F_2; t_2, r_2, \ell _2),$ a bijection $\phi$ from $F_1$ to
$F_2$ is called a {\it map isomorphism} if $\phi t_1=t_2\phi,$ $\phi r_1=r_2\phi$ and $\phi \ell _1=\ell _2\phi.$ In particular,
if $\MM _1=\MM _2=\MM ,$ then $\phi $ is called an {\em automorphism} of $\MM.$
The automorphisms of $\MM $ form a group $\Aut (\MM)$  which is called the {\em automorphism group} of the map $\MM$. By the definition of map isomorphism, we have $\Aut(\MM )=C_{S_F}(\Mon (\MM)),$  the centralizer of $\Mon(\MM)$ in
$S_F$. It follows from the transitivity of $\Mon (\MM )$ on $F$ that   $\Aut (\MM )$ acts semi-regularly on $F.$  If the action is
regular, we call the map $\MM $ {\it regular}. As a consequence of a result in permutation group theory
(see \cite[I.Theorem 6.5]{Hup}), we get that for a regular map $\MM ,$ the two associated permutation groups $\Aut (\MM )$ and $\Mon (\MM )$
can be viewed as the right   regular  representations $R(G)$ and left regular representations  $L(G)$ of an abstract group $G\cong\Aut(\MM )\cong\Mon (\MM )$ mutually centralizing each other in $S_F$.
If there is no confusion we also identify  the elements $R(g)$ and $L(g)$  with  the elements $g\in G$,
 so that $\MM\cong \MM(G; t, r, \ell )$, which  is called an {\it algebraic map}.
Thus the subgroups $\lg t, r\rg$, $\lg t, \ell \rg $ and $\lg r,\ell \rg $, respectively, stand for the
stabilizer of a vertex, an edge and a face which are mutually incident.
Finally, two  maps $\MM(G; t_1,r_1, \ell _1)\cong \MM(G; t_2, r_2, \ell _2)$ if and only if there
 exists an automorphism $\s $  of $G$ such that $t_1^\s =t_2,$ $r_1^\s=r_2$ and $\ell _1^\s =\ell _2$ (see \cite{Jon2}).

\vskip 3mm (2) {\bf Orientably-regular maps.}
\vskip 3mm
An  orientably-regular map can be described in the following direct way.
 \begin{defi}\label{map2}
{\rm For a given finite set $D$ and two fixed-point-free  permutations $r, \ell $ on $D$ where $\ell $ is an involutory,  a triple
$\MM=\MM(D; r, \ell )$ is called a {\it combinatorial orientable   map} if $\lg r, \ell \rg $ acts transitively  on $D$.}
\end{defi}
For a given map $\MM=\MM(D; r, \ell ),$  $D$  is called the {\it arc} set,
$r, \ell$ are called {\it local rotation and arc-revision involution,} respectively.
The group $\lg r, \ell \rg $ is called the {\it monodromy group} of $\MM$, denoted by $\Mon(\MM)$.
We define the {\it vertices, edges} and {\it face-boundaries} of $\MM$ to be the
orbits of the cyclic subgroups $\lg r\rg$, $\lg \ell \rg $ and $\lg r\ell \rg $, respectively.
The incidence in $\MM $ is defined by a nontrivial set  intersection.

Similarly, one may define isomorphisms and automorphisms of $\MM$ and know that
  $\Aut (\MM )$ acts semi-regularly on $D$.  If the action is
regular, we call the map $\MM $ {\it regular}. In this case,
  the two associated permutation groups $\Aut (\MM )$ and $\Mon (\MM )$
can be viewed as the right   regular  representations $R(G)$ and left regular representations  $L(G)$ of an abstract group $G=\lg r, \ell\rg \cong\Aut(\MM )\cong\Mon (\MM )$ mutually centralizing each other in $S_D$ so that the map is denoted by $\MM(G; r, \ell )$.
Moreover,  two  maps $\MM(G; r_1, \ell _1)\cong \MM(G; r_2, \ell _2)$ if and only if there
 exists an automorphism $\s $  of $G$ such that $r_1^\s=r_2$ and $\ell _1^\s =\ell _2$.

\subsection{Some Known Results }
\begin{prop}
\label{nc}
{\rm (\cite[Chap.1, Theorem 6.11]{SUZ})}
Let $H$ be a subgroup of a group $G$.
Then $C_G(H)$ is a normal subgroup of
$N_G(H)$ and the quotient $N_G(H)/C_G(H)$ is isomorphic
to a subgroup of $\Aut (H)$.
\end{prop}

\begin{prop} \label{nil}{\rm \cite[VI. Hauptsatz 4.3]{Hup}} \label{nil} Let $G=N_1N_2\cdots N_k,$ where
$N_i$ is a nilpotent subgroup of $G$ for all $i\in \{ 1, 2, \ldots, k\} ,$ and $N_iN_j=N_jN_i$
for any $i,$  $j.$ Then $G$ is solvable.
 \end{prop}

\begin{prop}\cite[Theorem 1]{Gur}\label{primep}
   Let T be a nonabelian simple group with a subgroup $H<T$ satisfying
 $|T:H|=p^a,$ for $p$ a prime. Then one of the following holds:
\begin{enumerate}\item[{\rm (i)}]  $T=A_n$ and $H=A_{n-1}$ with
$n=p^a;$
\item[{\rm (ii)}]  $T=\PSL(n,q)$, $H$ is the stabilizer
of a projective point or a hyperplane in $\PG(n-1,q)$ and
$|T:H|=(q^n-1)/(q-1)=p^a;$ \item[{\rm (iii)}]  $T=\PSL(2,11)$ and
$H=A_5;$ \item[{\rm (iv)}]  $T=M_{11}$ and $H=M_{10};$ \item[{\rm
(v)}]  $T=M_{23}$ and $H=M_{22};$ \item[{\rm (vi)}] $T=\PSU(4,2)$
 and $H$ is a subgroup of index 27.\end{enumerate}
\end{prop}

\begin{prop} \cite[Theorem 1]{GW} \label{dihedral} Let $G$ be a finite group with dihedral Sylow 2-subgroups. Let
$O(G)$ denote the maximal normal subgroup of odd order. Then $G/O(G)$  is
  isomorphic to either a subgroup of ${\rm P\Gamma L}(2,q)$ containing $\PSL(2, q)$ where  $q$ is odd,
  or $A_7$, or a Sylow $2$-subgroup of $G$.
 \end{prop}

\begin{prop} {\cite{Hup}} \label{sol}  Let $\{p_1, p_2, \cdots, p_l\}$ be the set of prime divisors of $|G|$. Suppose $G$ is solvable. Then
\begin{enumerate}\item[{\rm (i)}]
$F(G)=O_{p_1}(G)\times O_{p_2}(G)\times \cdots \times O_{p_l}(G)$;
\item[{\rm (ii)}] $C_G(F(G))\le F(G)$;
\item[{\rm (iii)}]  $\Phi(G)\le F(G)$.
\end{enumerate}
 \end{prop}

\begin{prop} {\cite[A corollary of Theorem 4.7A]{DM}} \label{affine}   Every solvable  primitive group is an affine group.
 \end{prop}

Recall that a finite $p$-group $G$ is a {\em special $p$-group}
if  either  $G$ is an elementary abelian $p$-group; or
$\Phi(G)=G'=Z(G)$ is an elementary abelian $p$-group.
\f Furthermore,   a special $p$-group $G$ is  called {\em extra-special} if $|Z(G)|=p.$
We remind that
a group $G$ is said to be {\em the central product of $A_1$, $A_2$,$\cdots$, $A_n$ related to the central subgroup $Z$} if
 $G=A_1A_2\cdots A_n$ such that for any $1\leq i<j\leq n,$ we have $A_i\cap A_j=Z$ and $[A_i, A_j]=1$, while
the central product is denoted  by $G=A_1\ast A_2\ast\cdots\ast A_n.$

A {\it maximal class} $p$-group of order $p^n$ is the group which has nilpotent class $n-1$.

\vskip 3mm  The following proposition gives the structure of extra-special $p$-groups.

 \begin{prop} {\rm(\cite{Hup})} \label{esp}
 Let $G$ be an extraspecial $p$-group. Then   $|G|=p^{2m+1}$ for some integer $m\ge 1$, and $G$ is one of the following two types:
$$G\cong\underbrace{ N\ast\ldots\ast N}_{m\, {\rm times}}\,\,or\,\,G\cong \underbrace{N\ast\ldots\ast N}_{m-1\, {\rm times}}\ast M,$$
where for $p=2$, $N\cong \B{D}_8$ and $M\cong \B{Q}_8$; for $p>2$,
$N$ and $M$  are nonabelian groups of order $p^3$ with  the respective exponent $p$ and $p^2.$
\end{prop}

\begin{prop}{\cite[Theorem 4.22]{SUZ}} \label{Spec} Let $G$ be a $p$-group. Assume that every characteristic abelian subgroup of $G$ is cyclic. Then $G$ is a central product of two subgroups $E$ and $S$ such that
\begin{enumerate}\item[{\rm (i)}]
$E=\{1\}$ or $E$ is an extra-special $p$-group where ${\rm exp}(E)=p$ for $p>2$; and
\item[{\rm (ii)}] $S$ is cyclic or $2$-group of maximal class.

\end{enumerate} \end{prop}

\section{Solvability  of $p$-Maps}
From now on let $\MM$ be an orientably-regular (resp. regular) $p$-map of order $p^k$ and valency $n\ge 2$, where $p$ is odd and $k\ge 1$. Let $X$ be the underlying graph of $\MM$ and  $V$  the set of vertices of $X$.  Set $G=\Aut^+(\MM)$ (resp. $\Aut(\MM)$) and take a Sylow $p$-subgroup $P.$ Pick up a vertex $v\in V$ such that  $H:=G_v=\lg r\rg $ (resp. $\lg r, t\rg$).

The following lemma is a well-known fact and here we give a proof for a sake of completeness.
\begin{lem}
  $P$ acts transitively on $V$ so that $G=PG_v=PH$.
\end{lem}
 \demo Set $|G|=p^nm$ where $(p,m)=1$. For any $v\in V$, we have $|G|=|G_v|p^k$.  Since $|P_v|=|P\cap G_v|\leq p^{n-k}$, then $|PG_v|=|P||G_v|/|P\cap G_v|\geq p^nm=|G|$, forcing that $PG_v=PH=G$.\qed

\vskip 3mm
 In this section we show  Theorem~\ref{main}.(1), that is the solvability  of  $\MM.$
 To do that, it suffices to prove  the following  two lemmas  formulated in the language of abstract groups.

\begin{lem} \label{cyclic}
  Let  $G$ be a group having a cyclic  subgroup $H$  of index a $p$-power. Then $G$  is solvable.
\end{lem}
\demo
   As mentioned before,   $PH=G=HP$ and both $P$ and $H$ are nilpotent. By Proposition~\ref{nil}  again, $G$ is solvable.\qed

\begin{lem} \label{SC}
  Let  $G$ be a group having a subgroup $H\cong \B{D}_{2n}$ of index a $p$-power. Then $G$  is solvable.
\end{lem}
\demo Now $G=PH$, where $H\cong \B{D}_{2n}$. In what follows,  we deal with two cases when $p=2$ and  $p$ is odd, separately.

 \vskip 3mm
 {\it Case 1: $p=2$.}
 \vskip 3mm
 Let $C$ be the maximal subgroup of odd order in   $H\cong \B{D}_{2n}$.  Then  $C$ is cyclic and $G=PC=CP$. By Proposition~\ref{nil},  $G$ is  solvable.

 \vskip 3mm
 {\it Case 2: $p$ is odd.}
 \vskip 3mm

For the contrary,  suppose $G$ is insolvable. Let $G$ be  a minimal   counterexample.
 Suppose that $G$ is nonabelian simple. In  Lemma~\ref{primep},  a classification of finite nonabelian simple groups containing   a subgroup $H$ of index $p^k$ is given. Checking  it, we know that
 the subgroup $H$  cannot be  a  dihedral group.  So in what follows, we assume that  $G$  is not  simple.
 Then $G$ contains  a normal subgroup $N$ such that $1\lneq N\lneq G$. We shall show  the solvability of $G$  by  showing that both $G/N$ and $N$ are solvable.

\vskip 3mm {\it Step 1:} The solvability  of $G/N$.

From $|G:H|=|G:HN||HN:H|=p^k$ for some integer $k$, we know that  $|G:HN|$ is a $p$-power. Since $|G/N:HN/N|=|G:HN|$,  it follows that  $G/N$ has a subgroup $HN/N$ of  index a $p$-power..

Since $HN/N\cong H/H\cap N$ and $H\cong \B{D}_{2n}$, we get that   $HN/N$  is isomorphic  either to $\B{Z}_2$ or to a dihedral group.
  If $HN/N\cong \B{Z}_{2}$, then $G/N$ is solvable by Proposition~\ref{nil}. Suppose that  $HN/N$ is dihedral.   Since  our group $G/N$ satisfies the conditions of the lemma and   $|G/N|<|G|$,
   we get   that $G/N$ is solvable,  taking into account that $G$ is assumed to be  a minimal   counterexample.

\vskip 3mm {\it Step 2:} The solvability  of $N$.

Note that $N$ contains a subgroup  $H\cap N$  such that   $|N:H\cap N|=|HN:H|=p^{s}$ for some nonnegative  integer $s$.
 First suppose  that $H\cap N$ is cyclic. Take $P_1\in \Syl_p(N)$ so that $N=P_1(H\cap N)=(H\cap N)P_1$. By Proposition~\ref{nil},  $N$ is solvable.
 Second suppose that $H\cap N$ is dihedral. Since $|N|<|G|$, it follows from the minimality of $|G|$ that $N$ is solvable.\qed
\vskip 3mm

\section{Normality of $p$-Maps  where $p\ge 5$}
In this section we prove the main part of Theorem~\ref{main}.(2), that is, the normality of orientably-regular (resp. regular) $p$-maps where $p\ge 5$. Again, we  formulate the result in the language of abstract groups.
 It will be proved in Lemma~\ref{ortien}  that every    orientably-regular $p$-map $\MM$ is normal  if $p$ is odd; and  in Lemmas~\ref{O_2'}-\ref{2-group}  that  every  regular $p$-map $\MM$ is normal  if $p\ge 5$.
For some terminologies and notations used below, the readers may see them in the end of Section 1.

\begin{lem} \label{ortien}
Suppose $G=\langle r,\ell\rangle$  such that $\ell$ is an involution,   $|r|$ is even and $|G:\lg r\rg |=p^k$,  where $p$ is odd.  Then $G$ contains  a normal Sylow $p$-subgroup.
\end{lem}
\demo
Let $H=\lg r\rg$ and $P\in \Syl_p(G)$ so that $G=PH$ and $|r|$ is even. By Lemma~\ref{sol}, $G$ is solvable.  For the contrary,  assume that $P$ is nonnormal in $G$. Then $O_p(G)<P$.
Let $\og=G/O_p(G)$. Then  $O_p(\og)=1$. Set $H_1$ be the $p'$-Hall subgroup of $H$.   Then $\og =\op \oh=\op\oh_1$.

   Since $O_p(\og)=1$ and $|G|=p^k|H|$,  the Fitting subgroup  $F(\og)$ of $\og$ is contained  in $\oh_1$.  Since  $\og$ is solvable and $\oh_1$  is cyclic,  it follows from Proposition~\ref{sol}(ii)  that $\oh_1\le C_{\og}(F(\og))\leq F(\og)$, which gives    $F(\og)=\oh_1.$ In particular,
   $\oh_1\lhd \og$, forcing  $C_{\og}(\oh_1)\unlhd \og$.

Since $\oh\leq C_{\og}(\oh_1)$,  $C_{\og}(\oh_1)$ contains a Sylow $2$-subgroup of $\og$. Since $C_{\og}(\oh_1)\lhd \og$,   it contains all Sylow $2$-subgroups of $\og $ and
 then  contains all involutions of $\og$.  Now we have  $\oh, \overline{\ell}\le C_{\og}(\oh_1),$ which implies  $C_{\og}(\oh_1)=\og$.

     Finally, from $C_{\og}(\oh_1)=\og$, we get $[\oh_1, \op]=\o1$,  which implies   $\op\unlhd \op \oh_1=\og$, that is $P\lhd G$, a contradiction.
        \qed

\begin{lem} \label{O_2'} Suppose $G=\langle t, r, \ell\rangle$  such that $t,r, \ell$ are  involutions, $t\ell =\ell t$ and  $\langle r, t\rangle\cong \B{D}_{2n}$   where $|G:\lg r, t\rg |=p^k$ for an odd prime $p$.
 Then   $G/O_{2^{'}}(G)$ is isomorphic  to  either $S_4$ or    a Sylow $2$-group of $G$. Moreover, if $p\ge5$, then $G/O_{2^{'}}(G)\not \cong S_4$.
\end{lem}
\demo Suppose $4\nmid |G|$.   Then  $|G|=2m$ for some odd integer $m$ and so $G$ has a unique subgroup of order $m$, that is $O_{2'}(G)$. Clearly,  $G/O_{2'}(G)\cong \B{Z}_2$,  as desired. Remark that
this group corresponds to a degenerated case for regular maps when $\ell=t$.

Suppose $4| |G|$. Then $G$ contains  a dihedral Sylow $2$-subgroup.  By Proposition~\ref{dihedral},  we get $G/O_{2^{'}}(G)$ isomorphic to $A_{7}$, a group $N$ such that $\PSL(2,q)\leq N\leq  \PGammaL(2,q)$ with $q$  being an odd prime power or a Sylow $2$-subgroup of $G$.

Since $G$ is solvable and $\PSL(2,q)$ is simple when $q\ge 5$,  the first case was excluded and the second case  happens only if $N=\PSL(2,3)\cong A_4$ or ${\rm P\Gamma L}(2,3)\cong S_4$. Note that $A_4$ cannot be generated by its involutions,  we get that  $G/O_{2^{'}}(G)$ is isomorphic to either  $S_4$ or    a Sylow $2$-subgroup of $G$. \qed

 Suppose that  $G/O_{2^{'}}(G)\cong S_4$ where $p\ge 5$.  From  $PO_{2'}(G)/O_{2'}(G)\lessapprox S_4$, we get that $P\le O_{2'}(G)$ and so $G/O_{2'}(G)=\op \oh=\oh $, but $S_4$ can not be the homomorphic image of a dihedral group, a contradiction. \qed

\begin{lem} \label{2-group} With the notations in Lemma~\ref{O_2'}, suppose that $G/O_{2^{'}}(G)$ is isomorphic to  a Sylow $2$-subgroup of $G$, where $p$ is odd. Then $G$ contains a normal Sylow $p$-subgroup.
\end{lem}
\demo  Under the hypothesis, we have  $G=L\rtimes S$ where $L$ is a characteristic subgroup of odd order of $G$    and  $S\in \Syl_2(G)$ which is either $\B{Z}_2$ or dihedral. Then $L=PH_1$, where $H_1$ is the $\{p,2\}'$-Hall subgroup of $H$ and clearly, $H_1$ is cyclic.   Since $G=L\rtimes S$ and $L\lhd G$, we get that  $\Syl_p(L)=\Syl_p(G)$ and so $O_p(L)=O_p(G)$.
Remind that  by Lemma~\ref{sol}, $G$ is solvable.

Assume that $H_1=1$.  Then  $P=L\lhd G$. So in what follows, we assume that $H_1\ne 1$.
 Let $\og=G/O_p(G)$. For the contrary, suppose that $P\ntriangleleft G$, that is $\op\ne \o1$. Then we  first derive the following two facts:
\vskip 3mm
{\it Fact 1:  $F(\ol)={\overline H_1.}$ }
\vskip 3mm
Form $\og=G/O_p(G)=\ol \oh,$ we get  $O_p(\ol )=1$, which in turn implies the Fitting subgroup $F(\ol )$ of $\og$  is contained in  $O_{p'}(\ol )$, that is the intersection of all $p'$-Hall subgroups of $\ol $.
 Since  $H_1$ is a $p'$-Hall subgroup of $L$ and also a $\{p,2\}'$-Hall subgroup of $H$,  we have $O_{p'}(\ol )\leq \oh_1$.
   Since $\og $ is solvable and $\oh_1$  is cyclic, we have  from Proposition~\ref{sol}.(ii) that $$O_{p'}(\ol)\le \oh_1\le C_{\ol }(O_{p'}(\ol))\le  C_{\ol}(F(\ol))\leq F(\ol)\le O_{p'}(\ol),$$
   which forces   $\oh_1=F(\ol)\lhd \og.$
\vskip 3mm
{\it Fact 2: There exists a $\oq\in \Syl_q(F(\ol))$ for some $q\ne p$  such that $2p\mid |\og:C_{\og}(\oq)|$ }.
\vskip 3mm
Since  $C_{\ol}(F(\ol))\le F(\ol)=\oh_1$ (Fact 1) and $\op \varsubsetneq \oh_1$,  we get  $[\op, F(\ol)]\ne \o1$.
 Therefore, there exists some $\oq\in \Syl_{q}(F(\ol))$  such that $[\op, \oq]\ne \o1$ where $q\ne p$. Since $\oq \char F(\ol)\char \ol\lhd \og$,
  we have $\oq\unlhd\og$ and thus $C_{\og}(\oq)\unlhd\og$.  Therefore,   $C_{\og}(\oq) $ contains one of Sylow $p$-subgroups (resp. Sylow 2-subgroups) if and only if  it contains all Sylow $p$-subgroups  (resp. all involutions).
   Since  in $\oh\cong \B{D}_{2n}$, their exists an involution, say $\olh \in \oh$ such that  $[\oq, \olh]\ne 1$. Now  $\op $ and $\lg \olh \rg $  are not contained in  $C_{\og}(\oq)$. Therefore, $C_{\og}(\oq)$   contains neither a Sylow $p$-subgroup of $G$ nor a Sylow
   2-subgroup of $G$, which means   $2p\di |\og:C_{\og}(\oq)|$.

\vskip 3mm
Come back to our proof. By Proposition~\ref{nc}, we have
$$N_{\og}(\oq)/C_{\og}(\oq)=\og/C_{\og}(\oq)\lessapprox \Aut(\oq).$$
 Since $\Aut(\oq)$ is cyclic, we get that $\og/C_{\og}(\oq)$ is cyclic. Let $\ot /C_{\og}(\oq)$ be the Sylow 2-subgroup of  $\og/C_{\og}(\oq)$. Since $2p\di |\og:C_{\og}(\oq)|$ by Fact 2,  $\ot\ne \o1, \og$.
 Since $\ot $ is a subgroup of $\og$ of odd index,   $\ot $ contains a Sylow $2$-subgroup of $\og $. Since $\ot\lhd \og $, it contains all Sylow $2$-subgroups of $\og$, and in particular it contains all involutions.
   Then $\og =\lg \olr, \olt, \oll\rg \le \ot$,  a contradiction.

   In summary,  we get $P\lhd G$. \qed

\section{Normallity of  3-Maps}
In this section, we deal with 3-maps.

\begin{lem}\label{3-map}
Suppose that $\MM$ is an   orientably-regular 3-map or  a regular $3$-map. Then either
\begin{enumerate}
   \item[{\rm (1)}] $\MM$ is normal; or
     \item[{\rm (2)}]  $\MM$ is nonorientable and  nonnormal regular.  Moreover, it is a regular covering of ${\cal C}(3,2)$,  whose  covering transformation group is a $3$-group.
         \end{enumerate}
\end{lem}
\demo  From Lemma~\ref{ortien}, every  orientably-regular  $3$-map is  normal.

Let $\MM$ be a regular 3-map and set $G=\Aut(\MM)$.   If  $G/O_{2'}(G)$ is  isomorphic to a Sylow 2-subgroup of $G$,  then $\MM$ is normal by Lemma~\ref{2-group}.
Following Lemma~\ref{O_2'}, we need to consider the case  that  $G/O_{2'}(G)\cong S_4$. Let $P\in \Syl_3(G)$.
The proof will be divided into the following two steps.
\vskip 3mm
{\it Step 1:  Show $G=PD_8$ and $\MM$ is nonnormal and nonorientable.}
\vskip 3mm
Since $G/O_{2'}(G)\cong S_4$, we get $H=H_1S$ where $S\cong D_8$ and $H_1$ is the $\{3,2\}'$-Hall subgroup of $H$, which is  cyclic. Let $L=O_{2'}(G)$. Then
$G=PH=LPS$. Since $(|LH_1/L|, 2\cdot 3)=1$, we have $H_1\le L$. Since $LO_3(G)$ is a normal subgroup of odd order of $G$, we get $O_3(G)\le L$.

If $P\lhd G$, then $P\leq O_{2'}(G)$, which forces  $G/O_{2^{'}}(G)\not\cong S_4$, reminding $p=3$. Therefore  $P$ is nonnormal in $G$, that is $O_{3}(G)<P$.

Suppose that $H_1=1$. Then   $O_{2'}(G)=O_3(G)$ and  $G/O_3(G)\cong S_4$  so that   $G=PD_8$.

Suppose  $H_1\ne 1$. Then there exists a prime divisor $q$  of $|G|$  where $q\ne 2, 3$. Set $\og=G/O_3(G)$.
 Then we first show the following two facts.

\vskip 3mm
{\it Fact 1:  $F(\ol)={\overline H_1.}$ }
\vskip 3mm
  The proof is completely the same as that in Fact 1 of Lemma~\ref{2-group}.

 \vskip 3mm
{\it Fact 2: For any  $\oq\in Syl_q(F(\ol))$ where $q\ne 2, 3$, we have $6\mid |\og:C_{\og}(\oq)|$ }.
\vskip 3mm
 By Fact 1, $F(\ol)={\overline H_1.}$   Let $\oq\in \Syl_q(F(\ol))$, where we may take $Q\le H_1\le H$. 
 
  Since $\oq \char F(\ol)\char \ol$,  we have $\oq\unlhd \og$, 
 forcing  $\oc:=C_{\og}(\oq)\unlhd \og .$ 
 Since $[\oq,\os]\neq \o1$, some element of $\os\cong D_8$ is not contained in $\oc$.  By the normality of $\oc $, we know that $\oc$ cannot contain any Sylow 2-subgroup of $\og$, that is $2\di |\og:\oc|$. But from $\oq\leq \oh$, we know that a Sylow $2$-subgroup $\oa$ of $\oc$ is isomorphic to $\B{Z}_4$.

 Suppose that $\oc$ contains a Sylow 3-subgroup of $\og$. Then $3\di |\oc\ol/\ol|$.
  Since  $\oc\ol/\ol \lhd \og/\ol\cong G/L\cong S_4$, it follows that  $\oc\ol/\ol$ contains a subgroup isomorphic to $A_4$,   contradicting to $\oa\cong \B{Z}_4$.
      Therefore,  $\oc$  does not contain any Sylow 3-subgroup of $\og$,  that is $3\di |\og:\oc|$.

\vskip 3mm
Finally, using the same arguments as  that in  last paragraph of  Lemma~\ref{2-group}, we may get  $P\lhd G$, a contradiction. Therefore, $H_1$ must be trivial. In particular, $\og=G/O_3(G)\cong S_4.$ Note that $\langle\olr\olt,\olt\oll\rangle$ is a normal subgroup of $\og$ containing a subgroup $\langle\olr\olt\rangle\cong \B{Z}_4,$ we must have $\og=\langle\olr\olt,\olt\oll\rangle$. On the other hand, that $|G:\langle rt,t\ell\rangle|\leq 2$ yields $O_3(G)\leq \langle rt,t\ell\rangle$, forcing $G=\langle rt,t\ell\rangle$. Hence $\MM$ is a nonorientable map.
\vskip 3mm
{\it Step 2:  Show ${\overline \MM}={\cal C}(3,2)$.}
\vskip 3mm
  Consider the quotient 3-map ${\overline \MM}$ induced by $O_3(G)$. Then $\og =\lg \olr, \olt, \oll\rg =S_4$ where $\oh=\lg \olr, \olt\rg =\B{D}_8$ and so ${\overline \MM}$ has three vertices. Since the kernel of the action of $\og$ on vertices  is isomorphic to $\B{D}_4$, there are double edges (four flags) between every two adjacent vertices  and thus the underlying graph is $C_3^{(2)}$ and ${\overline\MM}\cong {\cal C}(3,2)$. Therefore,  $\MM$ is a regular cover of ${\cal C}(3,2)$, with the covering transformation group  $O_3(G)$, a 3-group.\qed

\vskip 3mm
To end up this subsection, we give two  examples that $|O_3(G)|=1$ or 3.
\begin{exam}   $O_3(G)=1$:  Let $G=S_4$. Take $r=(13), \ t=(12)(34)$ and $\ell =(12)$. Then $\MM(G; r, t, \ell)$ is a nonorientable nonnormal 3-map which has three vertices and six edges. Since $|\langle r,\ell\rangle|=6$, the map has four faces. The genus  $g$ of the map is $1$.

\end{exam}

\begin{exam}  $|O_3(G)|=3$: Let $G=(\langle b\rangle\times \langle c\rangle)\rtimes \langle d,e\rangle$ with the defining relations
$$b^2=c^2=e^2=d^9=[b,c]=1,  b^d=c,c^d=bc,    d^e=d^{-1}, b^e=c, c^e=b.$$
\f Then $O_3(G)=\langle d^3\rangle$ and $ G/O_{3}(G)\cong S_4$. Let $r=e,t=b,\ell=de$, then $\MM(G; r, t, \ell)$ is a nonorientable nonnormal 3-map which has $9$ vertices, $18$ edges and $4$ faces. The genus  $g$ of the map is $7$.

  \end{exam}

\section{Normality of  2-Maps}
\begin{lem} \label{2-map} Suppose that $\MM$ is an   orientably-regular 2-map or a regular $2$-map. Then  either
\begin{enumerate}
   \item[{\rm (1)}] $\MM$ is normal; or
   \item[{\rm (2)}] $\MM$  is a regular covering of one  of  the following maps ${\overline \MM}$,  whose   covering transformation group is a $2$-group:
   \begin{enumerate}  \item[{\rm (2.1)}] ${\overline \MM}={\cal D}(m,e)$, where $m\ge 3$ is odd and $e^2\equiv 1(\mod m)$ but $e\not\equiv 1(\mod m)$.
     $\MM$ is   nonnormal orientably-regular , more precisely, it is either a chiral map  or reflexible and  nonnormal  regular;

      \item[{\rm (2.2)}] ${\overline \MM}=DM(m)$ and $\MM$ is nonorientable and nonnormal regular;

    \item[{\rm (2.3)}]  ${\overline \MM}=EM(m)$ and $\MM$ is normal  orientably-regular  but nonnormal regular.
     \end{enumerate}   \end{enumerate}
        \end{lem}
\demo Suppose that  $\MM$ is a nonnormal orientably-regular 2-map or a regular $2$-map. Set $G=\Aut(\MM)$ or $\Aut^+(\MM)$ and take $P\in \Syl_2(G)$. Then  $G$ cannot be a $2$-group and $O_2(G)\lvertneqq P$. We shall deal with two cases, separately.
\vskip 3mm
 (1)  $\MM$ is a  nonnormal  orientably-regular $2$-map.
\vskip 3mm
   Set  $G=\lg r, \ell\rg $ and  $\og =G/O_2(G)$. Then $\og =\op \lg \olr\rg $, and the Fitting subgroup $F(\og)$ is of odd order, which implies  $F(\og)\le  \lg \olr\rg$.
Since $\og $ is solvable,  we have  $\lg \olr\rg \le C_{\og}(F(\og))\le F(\og)$, that is  $F(\og)=\lg \olr\rg $, and in particular  $|\olr|$ is odd.  Therefore,
$\og=\lg \olr , \oll\rg =\lg \olr \rg \rtimes \lg \oll\rg $. Since $\op\ne \o1$,  we get $\oll\ne \o1$  and so  $|P:O_2(G)|=2$. This   implies  that the quotient map ${\overline \MM}$ of $\MM$ induced by $O_2(G)$  has two vertices
and then the underlying graph of ${\overline \MM}$  is  a dipoles ${\cal D}_m$, where $m=|\olr|$ is odd.
A classification of orientably-regular maps of ${\cal D}_m$ was  given in \cite{NS}. As mentioned in Section 1,  every such map is isomorphic to  ${\cal D}(m,e)$, where $e^2\equiv 1(\mod m)$.
  Since $\og$ is nonabelian, $e\not\equiv 1(\mod m)$.   Therefore, $\MM$ is a regular cover of ${\cal D}(m,e)$,   with the covering transformation group $O_2(G)$, a 2-group.

 \vskip 3mm
(2)  $\MM$ is a nonnormal  regular $2$-map.
\vskip 3mm
Set  $G=\lg r, t, \ell\rg $   and $H_1$  the $2'$-Hall subgroup of $\lg r, t\rg $.
 Let $\og =G/O_2(G)$. Then $\og =\op \lg \olr, \olt\rg $. Then $|F(\og)|$ is odd and so $F(\og)\le \oh_1$.   From $\lg \olr\olt\rg \le C_{\og}(F(\og))\le F(\og)$, that is  $F(\og)=\lg \olr\olt\rg \lhd \og$, and thus $m:=|\lg \olr\olt\rg|$ is odd and  $\lg \olr, \olt\rg \cong \B{D}_{2m}$. Therefore,
$\og=\lg \olr \olt\rg \rtimes \lg \olt, \oll\rg =\lg \olr, \olt \rg \lg \oll\rg $, reminding $[t,\ell]=1$.

 Suppose that $\oll\not\in \lg\olr, \olt\rg $. Then
  ${\overline \MM}$  has two vertices and so its underlying graph  is a dipoles ${\cal D}_m$ where
$m$ is odd.   If  $\MM$ is nonorientable,  then   ${\overline \MM}$ is nonorientable too. However, By  \cite[Lemma 3.1]{HNSW}, the  dipoles ${\cal D}_m$ has a nonorientable map if and only if $m=2$. Since $m$ is odd,
our  $\MM$ is  orientable and (so reflexible).  Since $\og \cong \B{Z}_m\rtimes \B{D}_4$ and $F(\og)=C_{\og}(F(\og))\cong \B{Z}_m$, $m$ must contains at least two prime divisors by Proposition 2.3. Clearly, $\overline{\lg rt, t\ell\rg}$ has index 2 in $\og$ and so  $O_{2}(G)\le \lg rt, t\ell\rg $.
  Therefore,  $\MM$ is nonnormal orientably-regular if and only if it is nonnormal regular.

 Suppose that  $\oll\in \lg\olr, \olt\rg $. Then $\og \cong \B{D}_{2m}$ and so the quotient map ${\overline \MM}$ induced by $O_2(G)$  has only one vertex. Since $\oll \olt=\olt\oll$ and $\og \cong \B{D}_{2m}$ where $m$ is odd, we have either $\oll=\o1$ or $\oll=\olt$:
  \vskip 3mm
 (1) $\oll=1$: Following  \cite{LS}, ${\overline \MM}$ is said to be degenerate. The underlying graph is a $m$-semistar and it can be embedded in a disc, which is an orientable surface with  boundary. This map is denoted by $DM(m)$. Since $\lg \olr\olt, \olt\oll\rg =\lg \olr, \olt\rg =\og$, we  have $\lg rt, t\ell\rg =G$ and so $\MM$ is nonorientable.

\vskip 3mm
 (2) $\oll=\olt$:  Following  \cite{LS}, ${\overline \MM}$ is said to be redundant.  The underlying graph is a $m$-semistar too  and it can be embedded in a sphere.  This map is denoted by $EM(m)$. Since $\lg \olr\olt, \olt\oll\rg =\lg \olr\olt\rg \ne \og$, we  have $\lg rt, t\ell\rg \ne G$ and $O_2(G)\le \lg rt, t\ell\rg $  and so $\MM$ is orientable and so reflexible. Moreover, $O_2(G)\le \Syl_2(\lg rt, t\ell\rangle)$, which  means
 that $\MM$ is  normal orientably-regular  but nonnormal regular.
\qed

\vskip 3mm

In what follows, we construct four $2$-maps, which give  respective four examples satisfying the conditions in  Lemma~\ref{2-map}.(2).
\begin{exam}
Let $$G=((\langle a\rangle\times \langle b\rangle)\rtimes (\langle c\rangle\times \langle d\rangle))\rtimes \langle e,f\rangle \cong ((\B{Z}_4\times \B{Z}_4)\rtimes (\B{Z}_2\times \B{Z}_2))\rtimes S_3,$$
whit the following defining relations:
$$\begin{array}{ll} &a^4=b^4=[a,b]=1; c^2=d^2=[c,d]=1, a^c=ab^2, b^c=a^2b^3, a^d=a^3, b^d=b^3; \\
 &e^2=f^3=1,  f^e=f^2,  a^e=a^2b^3, b^e=a^3b^2 , c^e=cd, d^e=d,
a^f=a^2b, b^f=ab ,\\ &c^f=c, d^f=d.\end{array}$$
\f Clearly, $O_2(G)=\lg a, b, c, d\rg$  which is of order $2^6$ and $\langle e,f\rangle\cong D_6$.

Set $r=a^{-1}b^{-1}cf$ and $\ell=e$. Omitting computations in details, we get
\begin{eqnarray}
\label{x1}
 (i))\, [r,\ell]=df, \, (ii)\, r^2=af^{-1},\, (iii)\, (r\ell)^2(\ell r)^2=a^{-1}b,\,(iv)\,  r^3=a^2b^2c.\end{eqnarray}
Then (i) gives  $d, f\in \lg r, \ell\rg $, and combining  (ii), (iii) and (iv) in order we have  $a, b, c\in \lg r, \ell\rg $. So $G=\lg r, \ell\rg $.
Moreover, we state   that there exists no  $\s\in \Aut(G)$ such that $r^\s=r^{-1}$ and $\ell^\s=\ell$. Assume the converse.
Then from Eq(\ref{x1}) again, we have
$$f^\s=([r,\ell]^4)^\s=[r^{-1}, \ell]^4=a^2bf^{-1},\,\, a^\s=r^{-2}f^\s= r^{-2}a^2bf^{-1}=b,\,\, b^\s=a\  ({\rm as}\, \s^2=1)$$
 \f
 From $$ (a^2b^2)c^\s=(a^2b^2)^\s c^\s =(r^3)^\s=r^3=a^2b^2c,$$
  we have $c^\s=c$.  Then we get the following contradiction.
  $$r^\s=(a^{-1}b^{-1}cf)^\s=a^{-1}b^{-1}cf^\s=a^3b^2cf^{-1}\quad {\rm but} \quad r^{-1}=a^3cf^{-1}.$$

The above arguments tell us that
$\MM=\MM(G;r,\ell)$ is a chiral and nonnormal  $2$-map  whose quotient map $\overline{\MM} \cong {\cal D}(3,2)$, which is one of the cases in  Lemma~\ref{2-map}.(2.1). Finally, one may
check that $|r|=6$ and $|r\ell|=4$. Then  $\MM$ has $2^6$ vertices, $3\cdot 2^5$ faces and $3\cdot 2^6$ edges, which implies $g=1-\frac 12(2^6+3\cdot 2^5-3\cdot 2^6)=17$.
\end{exam}

  \begin{exam}
 In $G=\GL(2,3)$, set
  $$r={\small \left(\begin{array}{cc}-1&1 \\ 0 &-1\\\end{array}\right)},  \, \,  \ell={\small \left(\begin{array}{cc} 0&1 \\ 1 &0\\\end{array}\right)},  \, \,
  z={\small \left(\begin{array}{cc} -1&0 \\ 0&-1\\\end{array}\right)},  \, \,  w={\small \left(\begin{array}{cc} 1&0 \\ 0 &-1\\\end{array}\right)}.$$
 Then $G=\lg r, \ell\rg$. If a relation  $c(r,\ell)=1$, then $\ell$ must appear even number times in the word. So  $c(r,\ell)=1$
  if and only if $c(r, z\ell)=1$,  it follows that $\s: r\to r, \ell\to z\ell$ can be extended to an automorpism of $G$.
  Moreover,  $\s \cdot \Inn(w) $ fixes $\ell$ and maps $r$ to $r^{-1}$.
  Check that $|r|=6$ and $|r\ell|=8$. Then  the  $2$-map $\MM =\MM(G; r,\ell)$ is  reflexible and   has $8$ vertices, $24$ edges and $6$ faces, which implies that $g=12$.
  Since $O_2(G)\cong Q_8$, we have $G/O_2(G)\cong \B{D}_6$ and $\langle r\rangle O_2(G)/O_2(G)\cong \B{Z}_3$.  Thus the quotient map ${\overline \MM}$  has exactly two vertices and three edges,
   which implies ${\overline \MM}\cong {\cal D}(3,2)$.
\end{exam}

\begin{exam}
Let $G=S_4=\langle r,t,\ell\rg $, where $r=(12)$, $t=(13)$ and $\ell=(13)(24)$ and $\MM=\MM(G; r,t ,\ell)$. Since $\lg rt, t\ell\rg =G$, we know $\MM$ is a nonorientable and  nonnormal regular $2$-map,
which has $4$ vertices, $3$ faces and $6$ edges. Therefore, it is an embedding of $K_4$ into the projective plane. Moreover, $O_2(G)\cong \B{D}_4$ and ${\overline \MM}=DM(6)$.
    \end{exam}

\begin{exam}
Let $G=S_4=\langle r,t,\ell\rg $, where $r=(12)$, $t=(13)$ and $\ell=(24)$ and $\MM=\MM(G; r,t ,\ell)$. Since $\lg rt, t\ell\rg =A_4$, we know $\MM$ is a nonnormal orientably-regular  $2$-map,
which has $4$ vertices, $4$ faces and $6$ edges. Therefore, it is an embedding of $K_4$ into the sphere.    Moreover, $O_2(G)\cong \B{D}_4$ and  ${\overline \MM}=EM(6)$.
    \end{exam}

\section{Affine $p$-Maps}
By Theorem~\ref{main}, every orientably-regular $p$-map or regular $p$-map is normal, with some exceptions for $p=\{ 2 ,3\}$. In this section, we concentrate on normal $p$-maps.

 Suppose $\MM$ is a normal orientably-regular $p$-map or a normal regular $p$-map and let $G=\Aut^+(\MM)$ or $\Aut(\MM)$. Take $P\in \Syl_p(G)$. Then
$P\lhd G$. In general case,  a $p$-group may contain some characteristic subgroups, which  induce blocks on the set $V$ of all vertices, in other words, the induced action of $G^V$ on $V$  is imprimitive.
 Therefore, it is a framework to determine normal $p$-maps such that $G^V$ is primitive.

We use $c(P)$ to denote the nilpotent class of $P$ and set $\Omega_1(P):=\langle x\in P| x^p=1\rangle$, note that $\Omega_1(P)=\{x\in P| x^p=1\}$ when $p>2$ and $c(P)\leq 2$. Recall that $H_G$ denotes the core of $H$ in $G$.
\begin{prop} \label{pri1}
Suppose that $\MM$ is a normal  orientably-regular $p$-map with $p^k$ vertices  where $k\geq 2$. Let $G=\Aut(\MM)=\lg r, \ell\rg $,  $H=\langle r\rangle$ a point stabilizer, $P\in \Syl_p(G)$ and $P_0=H_G\cap P$.
 Suppose that $G^V$ is  primitive.    Then  we have
 \begin{enumerate}
   \item[{\rm (1)}] $P=P_0\times T$, where $T\cong \B{Z}^k_p$; or
   \item[{\rm (2)}] $P=E\ast P_0$, where $P_0=Z(P)$ and  $E$  is an extra-special $p$-group, where $E$  can be chosen as $\Omega_1(P)$ if $p>2$. In particular, $k$ is always even in this case.
 \end{enumerate}
\end{prop}
\demo  Suppose that $G^V$ acts  primitively on $V$. Since $G$ is solvable,  $G^V$ is solvable too. By Proposition~\ref{affine}, $G^V$   must be of affine type and $V$ can be identified with a vector space of dimension $k$.  Thus $\lg \olr\rg$ acts irreducibly on $V$ , $|\olr|\di {p^k-1}$ and $|\olr|\nmid {p^s-1}$ for any $1\leq s<k$, which implies
  $P/P_0\cong \B{Z}^k_p$ so that $G/H_G=\op \rtimes \lg \olr\rg $. Moreover, $P_0\le Z(G)$, as  $H/H_G\le C_G(P_0)/H_G\lhd G/H_G\cong \B{Z}_p^k\rtimes H/H_G$, a primitive group.

   In what follows we deal with two cases separately.
\vskip 3mm
{\it (1) $P$ is abelian}
\vskip 3mm
Since $P$ is abelian,  $P_0$ is cyclic and $P/P_0\cong \B{Z}_p^k$, we know that either $P=P_0\times \B{Z}^k_p$ or  $P=P_1\times \B{Z}^{k-1}_p$ where $P_1$ is cyclic and  $|P_1:P_0|=p$.  In fact, the second case cannot happen. Otherwise,
  $\Omega_1(P)$ is a characteristic subgroup of $P$, which induces blocks. By the primitivity of $G^V$, we get that   either $\Omega_1(P)$ fixes $V$ pointwise which means $\Omega_1(P)\le P_0$ or  $\Omega_1(P)$ acts  transitively on $V$,  but these two cases cannot happen. So $P=P_0\times \B{Z}^k_p$, as shown in  (1) of the lemma.

\vskip 3mm
{\it (2) $P$ is nonabelian}
\vskip 3mm
 Since a characteristic abelian subgroup is not transitive on $V$ and it induces blocks on $V$, it must be contained in   $P_0$. In  particular,  $Z(P)\le P_0$, and then $Z(P)=P_0$, which is cyclic.  By Proposition~\ref{Spec}  we deduce that $P$ is a central product of two groups $S$ and $E$, where $S$ and $E$ are given below according to $p>2$ and $p=2$, separately.

\vskip 3mm

{\it Step 1: $p>2$.}
\vskip 3mm
By Proposition~\ref{Spec}, $S$ is a  cyclic subgroup and $E$ is an extra-special $p$-group of exponent $p$.

If $S=1$, then $P$ is an extra-special $p$-group, as desired.  Suppose that $S>1$. First note that $S\leq Z(P)=P_0$, we thus have $Z (P)=S(Z(P)\cap E)$ by Dedekind's lemma. But from $1<|Z(P)\cap E|\leq |Z(E)|=p$, we derive that $Z(P)\cap E\cong \B{Z}_p$, which implies that $S=Z(P)=P_0$. Since ${\rm exp}(E)=p$, we have $E\leq \Omega_1(P)$ by definition. If there are elements $e\in E, s\in S$ such that $|es|\leq p$, then we must have both $|e|\leq p$ and $|s|\leq p$. Since any element $s\in S$ with $|s|\leq p$ is contained in $E$, we conclude that any element of order $p$ in $P$ lies in $E$, thus $\Omega_1(P)=E$, as required.

\vskip 3mm
{\it Step 2:  $p=2$}
\vskip 3mm
By  Proposition~\ref{Spec}, $E$ is either identity or an extra-special $2$-group; and $S$ is either a cyclic $2$-group or a $2$-group of maximal class.

If $|P_0|=2$, then $P$ is an extra-special $2$-group and we can write $P=P_0\ast P$. So we suppose that $|P_0|\geq 4$ from now on. If $P$ is a cyclic $2$-group or a maximal class $2$-group, then $P$ would contain an index $2$ characteristic cyclic subgroup containing $P_0$ which is neither transitive on vertices nor contained in the core of $H$, a contradiction. Thus  $E$ is an extra-special $2$-group.

Now we show that $S$ can not be a maximal class $2$-group. Assume the contrary. Set  $P_0=\lg x\rg $, where $x=es$ for some $e\in E$ and $s\in S$. Then, we have $e\in Z(E), s\in Z(S)$ and either $|e|\geq |x|\geq 4$ or $|s|\geq |x|\geq 4$, but the center of both  an extra-special $2$-group and a maximal class $2$-group are of order $2$, a contradiction.
Thus $S$ is cyclic and clearly we have $S\leq P_0$.

 Finally, we show that $S=P_0$. Since $S\leq P_0$, we have $P_0=S(P_0\cap E)$. Since $P_0=Z(P)$, we must have $1<|P_0\cap E|\leq |Z(E)|=2$, which implies that $P_0\cap E\cong \mathbb{Z}_2$ and then $S=P_0$.

 Note that $P/P_0\cong E/Z(E)$ and the order of an extra-special $p$-group is always of the form $p^{2m+1}$ for some positive integer $m$, the statement $k$  is even follows, as required.
\qed

\begin{prop} \label{pri2}
Suppose that $\MM$ is a  normal  regular $p$-map with $p^k(k\geq 2)$ vertices  such that   $G=\Aut(\MM)$ acts  primitively  on vertices. Let $H$ be a point stabilizer, $P_0=P\cap H_G$. Then $k$ is even and either
\begin{enumerate}
 \item[{\rm (1)}]  $P=P_0\times \B{Z}^k_p$; or
  \item[{\rm (2)}]   $P=\Omega_1(P)\ast P_0$, where $P_0=Z(P)$ and  $\Omega_1(P)=\{x\in P| x^p=1\}$  is an extra-special $p$-group.
\end{enumerate}
\end{prop}
\demo
Note that if $\MM$ is a normal regular $2$-map such that $G$ acts primitively on vertices, then $G$ is a $2$-group with $H$ maximal which is of less interest, thus our assumption $k\geq 2$ excludes that possibility. Now $p$ must be odd and the proofs for (1) and (2)  are similar to those in Lemma~\ref{pri1}, reminding a slight difference here is that we need to invoke Proposition 2.3 to prove that $P_0=Z(P)$ in case $P$ is nonabelian. Here we only need to prove that $k$ is even.
 It suffices to consider the case $G$ acts faithfully and primitively on vertices. Since $G$ is solvable, Proposition~\ref{affine} can be applied to show that $G$ must be of affine type.

Write $G=T\rtimes G_{\a}$ where $P=T=\B{Z}_p^k$ and  $G_{\a}$ is a point (vector) stabilizer which is an irreducible subgroup of $\GL(k,p)$. Let $C$ be the  cyclic subgroup of $G_{\a}$ with index 2.

Suppose that $C$ fixes some nonzero proper subspace of $V$, we can choose some $V_0\neq 0$ fixed by $C$ with minimal dimension, then $V^t_0$ is also a subspace fixed by $C$ for some $t\in G_{\a}-C$. Thus $V_0\cap V^t_0=0$ by assumption on $V_0$ and irreducibility of $G_{\a}$ implies that $V=V_0\oplus V^t_0$, hence ${\rm dim}(V)=2 {\rm dim}(V_0)$ and then $k$ is even.

Now we may suppose that $C$ is irreducible on $V$. By \cite[page 196,theorem 3.10]{Hup}, $C$ is conjugate to a subgroup of some Singer subgroup. However, the normalizer of any such subgroup in $\GL(k,p)$ is isomorphic to the semilinear group ${\rm \Gamma L}(1,p^k)\cong {\rm GF}(p^k)^*\rtimes \langle \sigma\rangle$ where $\langle \sigma\rangle$ is a cyclic subgroup of order $k$. This group contains a dihedral subgroup only if $k$ is even.
\qed
\vskip 3mm
In what follows, we construct three $3$-maps,   showing that all the cases in Propositions 7.1 and 7.2 may be realized.
\begin{exam}
Let $$G=\langle a\rangle\times ((\langle b\rangle\times \langle c\rangle\times \langle d\rangle)\rtimes \langle e\rangle)\cong \B{Z}_3\times ((\B{Z}_3\times\B{Z}_3\times\B{Z}_3)\rtimes \B{Z}_{26}) $$
with the following defining relations:
$$\begin{array}{ll}a^3=b^3=c^3=d^3=e^{26}=1,a^b=a^c=a^d=a^e=a,b^c=b^d=b,c^d=c;\\
b^e=d^2,c^e=bd^2,d^e=cd.\end{array}$$
Clearly the Sylow $3$-subgroup $P=\langle a\rangle\times (\langle b\rangle\times \langle c\rangle\times \langle d\rangle)\lhd G$.  Let $r=ae$ and $\ell=ce^{13}$. Then check   that  $G=\langle r,\ell\rangle$ and
 there exists no automorphism $\s$ of $G$ such that $r^\s=r^{-1}$ and $\ell^\s=\ell$. Set $H=\langle r\rangle$ where $H_G=\langle a\rangle$. Then $\MM:=\MM(G; r,\ell)$  is a normal orientably-regular and chiral $3$-map with $27$ vertices such that $\Aut^+(\MM)=G$ acts primitively on vertices.
\end{exam}

\vskip 3mm
\begin{exam}
Let $$G=(\langle a\rangle\times (\langle b\rangle\times \langle c\rangle))\rtimes \langle d,e\rangle\cong (\B{Z}_3\times (\B{Z}_3\times\B{Z}_3))\rtimes \B{D}_8 $$
with the following defining relations:
$$\begin{array}{ll}&a^3=b^3=c^3=d^4=e^2=1,a^b=a^c=a^d=a,b^c=b,b^d=b^2c^2,c^d=b^2c;\\
&a^e=a^2,b^e=bc,c^e=c^2,d^e=d^3.\end{array}$$
Clearly the Sylow $3$-subgroup $P=\langle a\rangle\times (\langle b\rangle\times \langle c\rangle)\lhd G$.   Let $r=ade,t=e,\ell=b^2cd^2e$ so that  $G=\langle r,t,\ell\rangle$ and $|G:\langle rt,t\ell\rangle|=2$. Set $H=\langle r,t\rangle$ so that $H_G=\langle a\rangle$. Then  $\MM:=\MM(G; r,t,\ell)$  is a normal regular  $3$-map with $9$ vertices such that $\Aut(\MM)=G$ acts primitively on vertices and clearly $\MM$ is also orientably-regular.
\end{exam}
\vskip 3mm
\begin{exam}
Let $$G=((\langle a\rangle\times \langle b\rangle)\rtimes \langle c\rangle)\rtimes \langle d,e\rangle\cong ((\B{Z}_3\times\B{Z}_3)\rtimes \B{Z}_3)\rtimes \B{D}_8 $$
with the following defining relations:
$$\begin{array}{ll}&a^3=b^3=c^3=d^4=e^2=1,a^b=a^c=a^d=a,b^c=ab,b^d=b^2c^2,c^d=b^2c;\\
&a^e=a^2,b^e=b^2,c^e=b^2c,d^e=d^3.\end{array}$$
Clearly the Sylow $3$-subgroup $P=\langle a\rangle\ast (\langle b,c\rangle)\lhd G$.   Let $r=ade,t=e,\ell=bcd^2$ so that  $G=\langle r,t,\ell\rangle=\langle rt,t\ell\rangle $. Set $H=\langle r,t\rangle$ so that $H_G=\langle a\rangle$. Then  $\MM:=\MM(G; r,t,\ell)$ is a normal regular  and nonorientable  $3$-map with $9$ vertices such that $\Aut(\MM)=G$ acts primitively on vertices.
\end{exam}

\vskip 3mm
\f {\bf Acknowledgments:}
This work is partially  supported by the National Natural Science Foundation of China (12071312).


{\footnotesize\begin{thebibliography}{99}





\bibitem{CCDKNW} D.A. Catalano, M.D.E. Conder, S.F. Du, Y.S. Kwon, R. Nedela and S.E. Wilson,
                 Classification of  regular embeddings of $n-$dimensional cubes,
                 {\it J. Algeb. Combin.}, {\bf 33}(2011)(2), 215--238.

\bibitem{CDNS} M.D.E., Conder, S.,F. Du, R. Nedela, M. \v Skoviera,  Regular maps with nilpotent automorphism group, {\it J. Algebraic Combin.}, {\bf 44}(2016), 863-874.

\bibitem{DJKNS1} S.F. Du, G.A.~Jones, J.H. Kwak, R.~Nedela and M.~\v Skoviera,
Regular embeddings of $K_{n,n}$ where $n$ is a power of $2$, I: Metacyclic case,
{\it  European J. Combin.}, {\bf 28}(6)(2007), 1595-1609.

\bibitem{DJKNS2} S.F. Du, G.A. Jones, J.H. Kwak, R. Nedela, M. \v Skoviera,
 Regular embeddings of $K_{n,n}$ where $n$ is a power of $2$, II: Nonmetacyclic case,
 {\it  European J. Combin.}, {\bf 31}(2010), 1946-1956.



\bibitem{DKN1} S.F. Du, J.H. Kwak and R. Nedela, A Classification of regular embeddings of graphs of order a product of two primes,
               {\it J. Algeb. Combin.}, {\bf 19}(2004), 123-141.





\bibitem{DK2} S.F. Du and J.H. Kwak, Nonorientable regular embeddings of graphs of order $p^2$,
{\it Discrete Math.}, {\bf 310}(2010), 1743-1751.

\bibitem{DM} J.D. Dixon and B. Mortimer, {\it Permutation Groups}, Springer-Verlag (1996).




\bibitem{GNSS} A. Gardiner, R. Nedela, J. \v Sir\' a\v n, M. \v Skoviera, Characterization of graphs which underlie regular maps on closed surfaces,
               {\it J. London Math. Soc.}, {\bf 59}(1)(1999), 100-108.


\bibitem{Gur} R.M. Guralnick, Subgroups of prime power index in a simple group, {\it J. Algebra}, {\bf 81}(1983), 304-311.

\bibitem{GW}  D. Gorenstein and    J. H. Walter, The characterization of finite groups with
dihedral Sylow 2-subgroups I, {\it J. Algebra}, {\bf 2}(1965), 85-151.

\bibitem{Hup} B.~Huppert, {\it Endliche Gruppen I}, Springer, Berlin, 1979.

\bibitem{HNSW} K. Hu, R.~Nedela, M.~\v Skoviera, N. Wang, Regular embeddings of cycles with
multiple edges revisited, {\it Ars, Math. Contemp.}, {\bf 8}(2015), 177-194.


\bibitem{JaJo} L.D. James and G.A. Jones, Regular orientable imbeddings of complete graphs, {\it J. Combin. Theory Ser. B }, {\bf 39}(1985) 353-367.


\bibitem{Jon1} G.A. Jones, Regular embeddings of complete bipartite graphs: classification and enumeration, {\it London Math. Soc.}, {\bf 101}(2010), 427-453.

\bibitem{Jon2} G.A. Jones, Maps on surfaces and Galois groups, {\it Math. Slovaca}, {\bf 47}(1997), 1-33


\bibitem{KK3} J.~H.~Kwak and Y.~S.~Kwon, Classification of nonorientable regular embeddings of complete bipartite graphs,
              {\it J. Combin. Theory Ser. B}, {\bf 101}(2011), 191-205.




\bibitem{KN1} Y.S. Kwon and R. Nedela, Non-existence of nonorientable regular embedings of $n$-dimensional cubes,
              {\it Discrete Math.}, {\bf 307}(2007)(3-5),  511-516.


 \bibitem{LS} C.H. Le and Jozef ~\v Sir\'{a}\v n, Regular maps whose groups do not act faithfully
on vertices, edges, or faces Classification of  regular embeddings of complete multipartite graphs, {\it   Euro. J.  Combin.}, {\bf 26}(2005), 521-541.




\bibitem{MNS}  A. Malni\v c, R. Nedela, and M.~\v Skovier,  Regular maps with nilpotent automorphism groups,
{\it European J. Combin.}, {\bf  33}(2012), 1974-1986.


     \bibitem{NS} R.~Nedela, M.~\v Skoviera, Regular embeddings of canonical double coverings of graphs,
              {\it J. Combina. Theory, Ser. B}, {\bf 67}(1996)(1-3), 249-277.



\bibitem{SUZ} M. Suzuki, {\it Group theory II}, Springer-Verlag, New York Berlin Heidelberg Tokyo, 1985.

\bibitem{Wil} S.E. Wilson, Cantankerous maps and rotary embeddings of $K_n$, {\it J. Combin. Theory Ser. B}, {\bf 47}(1989), 262--273.





\bibitem{ZDM1} Y.H. Zhu, W.Q. Xu, S.F.Du and X.S. Ma, On the orientable regular embeddings of order prime-cube,
             {\it  Discrete Math.}, {\bf 339}(2016), 1140--1146.
\end{thebibliography}}
\end{document}